\documentclass{icmart}

\contact[m.sapir@vanderbilt.edu]{Department of Mathematics\\
Vanderbilt University\\ Nashville, TN 37240}

\newcommand{\la}{\langle}
\newcommand{\ra}{\rangle}

\newcommand{\vk}{{van Kampen}\ }

\newcommand{\area}{{\mathrm{a}}}

\newcommand{\Cay}{\mathrm{Cayley}}
\newcommand{\sss}{{\mathcal S}}

\renewcommand{\thesubsubsection}{\thesubsection.\Alph{subsubsection}}
\newcommand {\N}{\mathbb{N}} 

\newcommand{\ct}{ contiguity }

\newcommand {\Z}{\mathbb{Z}}            

\newcommand {\me}{\medskip}

\newcommand {\iv}{^{-1}}

\newtheorem{theorem}{Theorem}[section]

\theoremstyle{definition}

\newtheorem{remark}[theorem]{Remark}

\title[Algorithmic and asymptotic properties of groups]{Algorithmic and
asymptotic properties of groups}

\author[Mark Sapir]{Mark Sapir\thanks{This work was supported in part by the NSF grants DMS
0245600 and DMS-0455881.}}

\begin{document}


\begin{abstract}
This is a survey of the recent work in algorithmic and asymptotic
properties of groups. I discuss Dehn functions of groups, complexity
of the word problem, Higman embeddings, and constructions of
finitely presented groups with extreme properties (monsters).
\end{abstract}



\section{Introduction}

Although the theory of infinite groups is very rich and full of
powerful results, there are very few results having more influence
on group theory and surrounding areas of mathematics (especially
geometry and topology) as the following five:

\begin{itemize}
\item The Boone-Novikov theorem about existence of finitely presented
groups with undecidable word problem \cite{Bo},\cite{No};

\item The Higman theorem about embeddability of recursively presented groups into
finitely presented groups \cite{Hi};

\item The Adian-Novikov solution of the Burnside problem \cite{AN};

\item Gromov's theorem about groups with polynomial
growth \cite{G2};

\item Olshanskii and his students' theorems about existence of groups with all proper
subgroups cyclic (Tarski monsters), and other finitely generated
groups with extreme properties \cite{book}.
\end{itemize}

In this paper, I am going to survey the last ten years of my work on
the topics related to these results.

\bigskip

{\bf Acknowledgement.} Most of the work surveyed here is joint with
J.-C. Birget, C. Dru\c tu, V. Guba, A. Olshanskii and E. Rips. I am
very grateful to them for cooperation.

\section{S-machines}

Recall that a Turing machine, say, with one tape is a triple $(Y, Q,
\Theta)$ where $Y$ is a tape alphabet, $Q$ is the set of states,
$\Theta$ is a set of commands (transitions) of the form
$\theta=[U\to V]$ where $U$ has the form $vqu$ and $V$ has the form
$v'q'u'$. Here $u,v,u',v'$ are words in the tape alphabet, $q,q'\in
Q$. A {\em configuration} of the Turing machine is a word $wqw'$
where $w, w'$ are words in the tape letters, $q$ is a state letter.
To apply the command $[U\to V]$ to a configuration, one has to
replace $U$ by $V$.

To specify a Turing machine with many tapes, one needs several
disjoint sets of state letters. A {\em configuration} of the machine
is a word of the form $u_1q_1u_2...u_Nq_Nu_{N+1}$ where $q_i$ are
state letters, $u_i$ are words in tape letters. Of course one needs
to separate tapes. That can be done by using the special symbols
(endmarkers) marking the beginning and the end of each tape. But
these symbols can be treated as state letters as well. Every
transition has the form $[U_1\to V_1,...,U_N\to V_N]$, where
$[U_i\to V_i]$ is a transition of a $1$-tape machine.

Among all configurations of a machine $M$, one chooses one {\em
accept} configuration $W$. Then a configuration $W_1$ is called {\em
accepted} if there exists a {\em computation} $W_1\to W_2\to...\to
W_n=W$ where each step consists in application of a command of $W$.

Recall that the time function of a (non-deterministic) Turing
machine is the smallest function $f(n)$ such that every accepted
input $w$ of size at most $n$ requires at most $f(n)$ steps of the
machine to be accepted.

The ``common denominator" of the proofs of most of the results I am
reviewing here is the notion of an $S$-machine that I introduced in
\cite{SBR}. Roughly speaking, $S$-machines make building groups with
prescribed properties as easy as programming a Turing machine.

Essentially, an $S$-machine is simply an HNN-extension of a free
group, although not every HNN-extension of a free group is an
$S$-machine.

Let us start with an example that we shall call the {\em Miller
machine}. It is the famous group of C. Miller \cite{Mi}. Let $G=\la
X\mid R\ra$ be a finitely presented group. The Miller machine is the
group $M(G)$ generated by $X\cup \{q\}\cup \{\theta_x\mid x\in
X\}\cup\{\theta_r\mid r\in R\}$ subject to the following relations
$$\theta x=x\theta,\quad \theta_x xq=qx\theta_x,\quad \theta_r q=qr\theta_r$$
where $\theta$ is any letter in $\Theta=\{\theta_x\mid x\in
X\}\cup\{\theta_r\mid r\in R\}$. Clearly, this is an HNN-extension
of the free group $\la X, q\ra$ with free letters $\theta\in
\Theta$. The main feature of $M(G)$ discovered by Miller is that
{\em $M(G)$ has undecidable conjugacy problem provided $G$ has
undecidable word problem.} In fact it is easy to see that $qw$ is
conjugated to $q$ in $M(G)$ if and only if $w=1$ in $G$.

To see that $M(G)$ can be viewed as a machine, consider any word
$uqv$ where $u, v$ are words in $X\cup X\iv$. If we conjugate $uqv$
by $\theta_r$, we get the word $uqrv$ because $\theta_r
q=qr\theta_r$ and $\theta_r$ commutes with $u$ and $v$ (here and
below we do not distinguish words that are freely equal). Hence
conjugation by $\theta_r$ amounts to executing a command $[q\to
qr]$. Similarly, conjugation by $\theta_x$ amounts to executing a
command $[q\to x\iv qx]$. If $u$ ends with $x$, then executing this
command means moving $q$ one letter to the left. Thus conjugating
words of the form $uqv$ by $\theta$'s and their inverses, we can
move the ``head" $q$ to the left and to the right, and insert
relations from $R$.

The work of the Miller machine $M(G)$ can be drawn in the form of a
diagram (see Figure \ref{fig1}) that we call a {\em trapezium}. It
is a tessellation of a disc. Each cell corresponds to one of the
relations of the group. The bottom layer of cells on Figure
\ref{fig1} corresponds to the conjugation by $\theta_x$, the next
layer corresponds to the conjugation by $\theta_r$, etc. These
layers are the so-called {\em $\theta$-bands}. The bottom side of
the boundary of the trapezium is labeled by the first word in the
computation ($uqv$), the top side is labeled by the last word in the
computation ($q$), the left and the right sides are labeled by the
{\em history of computation}, the sequence of $\theta$'s and their
inverses corresponding to the commands used in the {\em computation}
$uqv\to...\to q$. The words written on the top and bottom sizes of
the $\theta$-bands are the intermediate words in the computation. We
shall always assume that they are freely reduced.

\begin{center}
\begin{figure}[!ht]
\unitlength .7mm 
\linethickness{0.4pt}
\ifx\plotpoint\undefined\newsavebox{\plotpoint}\fi 
\begin{picture}(122.75,59)(-30,0)
\put(55.25,10.75){\framebox(19.75,8.5)[cc]{}}
\put(55.25,19.25){\framebox(19.75,8.5)[cc]{}}
\put(75,10.75){\framebox(7.25,8.5)[cc]{}}
\put(75,19.25){\framebox(7.25,8.5)[cc]{}}
\put(82.25,10.75){\framebox(7.25,8.5)[cc]{}}
\put(82.25,19.25){\framebox(7.25,8.5)[cc]{}}
\put(89.5,10.75){\framebox(7.25,8.5)[cc]{}}
\put(89.5,19.25){\framebox(7.25,8.5)[cc]{}}
\put(111,11){\framebox(7.25,8.5)[cc]{}}
\put(111,19.5){\framebox(7.25,8.5)[cc]{}}
\put(48,10.75){\framebox(7.25,8.5)[cc]{}}
\put(48,19.25){\framebox(7.25,8.5)[cc]{}}
\put(16,10.75){\framebox(7.25,8.5)[cc]{}}
\put(16,19.25){\framebox(7.25,8.5)[cc]{}}
\put(32.25,10.75){\makebox(0,0)[cc]{$\ldots$}}
\put(104.25,11.75){\makebox(0,0)[cc]{$\ldots$}}
\put(64.5,6){\makebox(0,0)[cc]{$xq$}}
\put(63.25,15.5){\makebox(0,0)[cc]{$qx$}}
\put(20,8){\makebox(0,0)[cc]{$y$}}
\put(19.5,15.5){\makebox(0,0)[cc]{$y$}}
\put(114.25,16){\makebox(0,0)[cc]{$z$}}
\put(115,8.5){\makebox(0,0)[cc]{$z$}}
\put(12,14.75){\makebox(0,0)[cc]{$\theta_x$}}
\put(122.75,15.5){\makebox(0,0)[cc]{$\theta_x$}}
\put(64.38,10.75){\vector(1,0){.07}}\put(62,10.75){\line(1,0){4.75}}
\put(64.63,19.25){\vector(1,0){.07}}\put(62.75,19.25){\line(1,0){3.75}}
\put(16,15.25){\vector(0,-1){.07}}\put(16,17.5){\line(0,-1){4.5}}
\put(118.13,15.63){\vector(0,-1){.07}}\multiput(118.25,17.75)(-.03125,-.53125){8}{\line(0,-1){.53125}}
\put(65.5,23.5){\vector(-2,-3){.07}}\multiput(68.5,27.75)(-.03370787,-.04775281){178}{\line(0,-1){.04775281}}
\put(11.25,24){\makebox(0,0)[cc]{$\theta_r$}}
\put(121,25){\makebox(0,0)[cc]{$\theta_r$}}
\put(19.5,25.75){\makebox(0,0)[cc]{$y$}}
\put(72,25.75){\makebox(0,0)[cc]{$x$}}
\put(113.75,26.25){\makebox(0,0)[cc]{$z$}}
\put(59.75,25.25){\makebox(0,0)[cc]{$qr$}}
\put(66,34.75){\makebox(0,0)[cc]{$\ldots$}}
\put(65.5,43.25){\makebox(0,0)[cc]{$\ldots$}}
\put(67,55.25){\vector(1,0){.07}}\put(59.5,55.25){\line(1,0){15}}
\put(66.25,59){\makebox(0,0)[cc]{$q$}}
\put(38.75,39.75){\vector(3,2){.07}}\qbezier(16,27.75)(39.75,38)(59.5,55.25)
\put(95.25,39.75){\vector(-3,2){.07}}\qbezier(118,27.75)(94.25,38)(74.5,55.25)
\end{picture}
\caption{Trapezium of the Miller machine for a deduction
$uqv\to...\to q$. Here $u=y...x$, $v=...z$.} \label{fig1}
\end{figure}
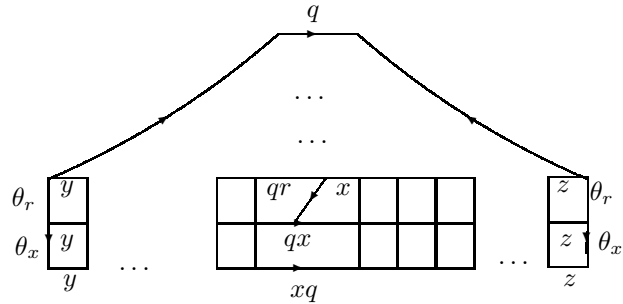
\end{center}

The Miller machine has one tape and one state letter. General
$S$-machines can have many tapes and many state letters. Here is a
formal definition.

Let $F(Q,Y)$ be the free group generated by two sets of letters
$Q=\cup_{i=1}^N Q_i$ and $Y=\cup_{i=1}^{N-1} Y_i$ where $Q_i$ are
disjoint and non-empty (below we always assume that $Q_{N+1}=Q_1$,
and $Y_N=Y_0=\emptyset$).


The set $Q$ is called the set of $q$-letters, the set $Y$ is called
the set of $a$-letters.

In order to define an HNN-extension, we consider also a collection
$\Theta$ of $N$-tuples of $\theta$-letters. Elements of $\Theta$ are
called {\em rules}. The components of $\theta$ are called {\em
brothers} $\theta_1,...,\theta_N$. We always assume that all
brothers are different. We set $\theta_{N+1}=\theta_1$,
$Y_0=Y_N=\emptyset$.

With every $\theta\in \Theta$, we associate two sequences of
elements in $F(Q\cup Y)$: $B(\theta)=[U_1,...,U_N]$,
$T(\theta)=[V_1,...,V_N]$, and a subset $Y(\theta)=\cup Y_i(\theta)$
of $Y$, where $Y_i(\theta)\subseteq Y_i$.

The words $U_i, V_i$ satisfy the following restriction:

\begin{itemize}
\item[(*)] For every $i=1,...,N$, the words $U_i$ and $V_i$ have the form
$$U_i=v_{i-1}k_iu_i, \quad V_i=v_{i-1}'k_{i}'u_{i}'$$ where $k_{i}, k_{i}'\in
Q_{i}$, $u_{i}$ and $u_{i}'$ are words in the alphabet $Y_{i}^{\pm
1}$, $v_{i-1}$ and $v_{i-1}'$ are words in the alphabet
$Y_{i-1}^{\pm 1}$.
\end{itemize}

Now we are ready to define an $S$-machine $\sss$ by generators and
relations. The generating set $X$ of the $S$-machine $\sss$ consists
of all $q$-, $a$- and $\theta$-letters. The relations are:

$$U_i\theta_{i+1}=\theta_i V_i,\,\,\,\, i=1,...,s, \qquad \theta_j a=a\theta_j$$
for all $a\in Y_j(\theta)$. The first type of relations will be
called $(q,\theta)$-{\em relations}, the second type -
$(a,\theta)$-{\em relations}.

Sometimes we will denote the rule $\theta$ by $[U_1\to
V_1,...,U_N\to V_N]$. This notation contains all the necessary
information about the rule except for the sets $Y_i(\theta)$. In
most cases it will be clear what these sets are: they are usually
equal to either $Y_i$ or $\emptyset$. By default $Y_i(\theta)=Y_i$.

Every $S$-rule $\theta=[U_1\to V_1,...,U_s\to V_s]$ has an inverse
$\theta\iv=[V_1\to U_1,...,V_s\to U_s]$; we set
$Y_i(\theta\iv)=Y_i(\theta)$.

\begin{remark} {\rm Every $S$-machine is indeed an HNN-extension of
the free group $F(Y,Q)$ with finitely generated associated
subgroups. The free letters are $\theta_1$ for every $\theta\in
\Theta$. We leave it as an exercise to find the associated
subgroups.}
\end{remark}

Every Turing machine $T$ can be considered as an $S$-machine $S'(T)$
in the natural way: the generators of the free group are all tape
letters and all state letters. The commands of the Turing machine
are interpreted as rules of the $S$-machine. The main problem in
that conversion is the following: there is a much bigger freedom in
applying $S$-rules than in executing the corresponding commands of
the Turing machine. Indeed, the Turing machine is in general not
{\em symmetric} (i.e. if $[U\to V]$ is a command of the Turing
machine then $[V\to U]$ is usually not) while every $S$-machine is
symmetric. Another difference is that Turing machines work only with
positive words, and $S$-machines work with arbitrary group words.
Hence the language accepted by $S'(T)$ is usually much bigger than
the language accepted by $T$.

Nevertheless, it can be proved that if $T$ is symmetric, and a
computation $w_1\to w_2\to...$ of the $S$-machine $S'(T)$ involves
only positive words, then that is a computation of $T$.

This leads to the following idea of converting any Turing machine
$T$ into an $S$-machine $S(T)$. First we construct a symmetric
Turing machine $T'$ that is equivalent to $T$ (recognizes the same
language). That is a fairly standard Computer Science trick (see
\cite{SBR}): the machine $T'$ first guesses a computation of $T$,
then executes it, then erases all the tapes. Note that the time
function and the space function of $T'$ are equivalent to the time
function of $T$.

The second step is to compose the $S$-machine $S'(T')$ with a
machine that checks positivity of a word. That machine starts
working after every step of $S'(T')$. That is if an application of a
rule of $S'(T')$ gives a non-positive (reduced) word then the
checking machine does not allow the machine $S'(T')$ to proceed to
the next step.

There are several checking machines. One of them - the {\em adding
machine} - is very simple but its time function is exponential (see
\cite{OSsmalldehn}). Another one is very complicated but it has a
quadratic time function (see \cite{SBR}).

Here is the definition of the adding machine. We present it here
also in order to show an example of a program of an $S$-machine. It
is not difficult to program an $S$-machine, but it does require some
practice.

Let $A$ be a finite set of letters. Let the set $A_1$ be a copy of
$A$. It will be convenient to denote $A$ by $A_0$. For every letter
$a_0\in A_0$, $a_1$ denotes its copy in $A_1$. The set of state
letters of the adding machine $Z(A)$ is $P_1\cup P_2\cup P_3$ where
$P_1=\{L\}, P_2=\{p(1),p(2),p(3)\}, P_3=\{R\}$. The set of tape
letters is $Y_1\cup Y_2$ where $Y_1=A_0\cup A_1$ and $Y_2=A_0$.

The adding machine $Z(A)$ has the following rules (there $a$ is an
arbitrary letter from $A$) and their inverses. The comments explain
the meanings of these rules.

\begin{itemize}

\item $r_1(a)=[L\to L, p(1)\to a_1\iv p(1)a_0, R\to R]$.
\me

{\em Comment.} The state letter $p(1)$ moves left searching for a
letter from $A_0$ and replacing letters from $A_1$ by their copies
in $A_0$.

\me

\item $r_{12}(a)=[L\to L, p(1)\to a_0\iv a_1p(2), R\to R]$.

\me

{\em Comment.} When the first letter $a_0$ of $A_0$ is found, it is
replaced by $a_1$, and $p$ turns into $p(2)$.

\me

\item $r_2(a)=[L\to L, p(2)\to a_0p(2)a_0\iv, R\to R]$.

\me

{\em Comment.} The state letter $p(2)$ moves toward $R$.

\me

\item $r_{21}=[L\to L, p(2)\to p(1), R\to R]$, $Y_1(r_{21})=Y_1,
Y_2(r_{21})=\emptyset$.

\me

{\em Comment.} $p(2)$ and $R$ meet, the cycle starts again.

\me

\item $r_{13}=[L\to L, p(1)\to p(3), R\to R]$, $Y_1(r_{13})=\emptyset,
Y_2(r_{13})=A_0$.

\me

{\em Comment.} If $p(1)$ never finds a letter from $A_0$, the cycle
ends, $p(1)$ turns into $p(3)$; $p$ and $L$ must stay next to each
other in order for this rule to be executable.

\item $r_{3}(a)=[L\to L, p(3)\to a_0p(3)a_0\iv, R\to R]$,
$Y_1(r_3(a))=Y_2(r_3(a))=A_0$

 \me

{\em Comment.} The letter $p(3)$ returns to $R$.

\end{itemize}

The underlying algorithm of the adding machine is simple: the
machine starts with a word $Lwp(1)R$, where $w$ is a word in $A\cup
A\iv$, $L, p(1), R$ are state letters. It considers the sequence of
indexes of the letters in $w$ as a binary number. The initial number
is $0$. The machine proceeds by adding 1 to this number until it
produces $2^n-1$ where $n$ is the length of the word (each cycle of
the machine adds a 1). After that, the machine returns the word to
its initial state (all indexes are 0). If the initial word contained
a negative letter, the state letter of the adding machine never
becomes $p(3)$.

To compose a checking machine $Z$ with an $S$-machine $\sss$ means
inserting state letters of $Z$ between any two consecutive state
letters of $\sss$, and changing the rules of $\sss$ in an
appropriate way: every rule of $\sss$ ``turns on" the checking
machines. After they finish their work, $\sss$ can apply another
rule (provided the word is still positive). If $Z$ is a checking
machine then the composition of $\sss$ and $Z$ is denoted by
$\sss\circ Z$.

The following results from \cite{SBR} are very important for the
applications. The equivalence of $S$-machines, their time functions,
space functions, etc. are defined as for ordinary Turing machines.

We say that two increasing functions $f,g\colon \Z_+\to \Z_+$ are
equivalent if

\begin{equation}\label{eq1} \frac1Cg\left(\frac{n}{C}\right)-Cn\le f(n)\le
Cg(Cn)+Cn\end{equation} for some constant $C$. We are not going to
distinguish equivalent functions in this note. Thus $n^{3.2}$ is the
same as $5n^{3.2}$ but different from $n^{3.2}\log n$.

\begin{theorem} (Sapir, \cite{SBR})) Let $T$ be a Turing machine. Then there exists an
$S$-machine $\sss$ that is polynomially equivalent to $T$. Moreover
the time function of $\sss$ is equivalent (in the sense of
(\ref{eq1})) to the cube of the time function of $T$, the space
function of $\sss$ is equivalent to the time function of $T$.
\end{theorem}

Moreover, one can use Miller's machines to simulate any Turing
machine.

\begin{theorem} (Sapir, \cite{OStalk}) For every Turing machine $T$
there exists a finitely presented group $G$ such that the Miller
machine $M(G)$ is polynomially equivalent to $T$.
\end{theorem}

Thus any Turing machine can be effectively simulated by an
$S$-machine with one tape and only one state letter.

\section{Dehn functions and the word problem}

\subsection{The definition} Let $G=\la X\mid R\ra$ be a finitely presented group. We
shall always assume that $X=X\iv$, $R$ is a collection of words in
the alphabet $X$ closed under taking inverses and cyclic shifts,
i.e. if $r\in R$, $r\equiv ab$ then $r\iv\in R$ and $ba\in R$.

The word problem in $G$ asks, given a word $w$ in $X$ (i.e. a
product of generators of $G$), whether $w$ is equal to $1$ in $G$.
Clearly, the word problem in all ``ordinary" groups is
algorithmically decidable. For example, if the group is linear, then
in order to check if $w=1$, one can just multiply matrices
representing the generators of $G$ in the order of their appearance
in $w$.

About 70 years ago, van Kampen noticed that $w=1$ in $G$ if and only
if one can tessellate a disc with boundary labeled by $w$ by tiles
(cells) whose boundaries are labeled by words in $R$.

That tessellation is called a {\em van Kampen} diagram for $w$. It
is also sometimes called {\em Dehn} or {\em disc} diagram of $w$.
For example, the trapezium on Figure \ref{fig1} is a van Kampen
diagram over the $S$-machine $M(G)$ with boundary label $huqvh\iv
q\iv$ where $h$ is the history of the computation.

For every $w=1$ in $G$, the area $\area(w)$ is the smallest number
of cells in the van Kampen diagram for $w$, or the simplicial area
of the (null-homotopic) loop labeled by $w$ in the Cayley complex
$\Cay(G,X)$. Combinatorially, that is the smallest number of factors
in any representation of $w$ as a product of conjugates of the words
from $R$. From the logic point of view, that is the length of the
shortest ``proof" that $w=1$ in $G$ (steps of the ``proof" are
insertions of relations from $R$ into $W$).

It is easy to see that the word problem in $G$ is decidable if and
only if the area of a word $w$ representing 1 in $G$ is bounded from
above by a recursive function in the length of $w$. Madlener and
Otto \cite{MO} and, independently, Gersten \cite{Ger} introduced a
very basic characteristic of the algorithmic complexity of a group
$G$, the Dehn function $\delta_G(n)$ of $G$: {\em it is the smallest
function $d(n)$ such that the area of a word $w$ of length $\le n$
representing 1 in $G$ does not exceed $d(n)$}. Of course
$\delta_G(n)$ depends on the choice of generating set $X$. But Dehn
functions corresponding to different generating sets are {\em
equivalent} in the sense of (\ref{eq1}). Similarly, one can
introduce the {\em isodiametric function} of $G$ by looking at the
diameter of \vk diagrams instead of their areas.

For example, the area of the trapezium on Figure \ref{fig1} is
approximately $|h|$ times the length of the longest $\theta$-band in
that trapezium, that can be interpreted as the product of the time
of the computation by its space. That observation is the key to
converting properties of the $S$-machine into the properties of the
Dehn function.

\subsection{The description}\label{td} Dehn functions reflect in an easy and natural way both
geometric and algorithmic properties of a group, so it is natural to
ask which functions appear as Dehn functions of groups.

The first observation is not difficult.

\begin{theorem}\label{2} (See \cite[Theorem 1.1]{SBR}) Every Dehn function of
a finitely presented group $G$ is (equivalent to) the time function
of a Turing machine solving (non-deterministically) the word problem
in $G$.
\end{theorem}

The proof of this  theorem in \cite{SBR} is more complicated than it
should have been. An easier proof can be obtained by using
\cite[Lemma 1]{OStalk}.

Not every increasing function can be equivalent to a time function
of a Turing machine. For example, if a time function $f(n)$ does not
exceed a recursive function then it must be recursive. On the other
hand, any ``natural" function is the time function of a Turing
machine. In particular, if $f(n)$ can be computed in time $\le f(n)$
then $f(n)$ is the time function of a deterministic Turing machine
computing $f(n)$.

By a theorem proved by Gromov and Olshanskii among others, every
finitely presented group with subquadratic Dehn function is in fact
hyperbolic, so its Dehn function is linear. It is possible to deduce
from a result of Kapovich and Kleiner \cite{KK} that if the Dehn
function is subquadratic even on an infinite subset of natural
numbers then the group is still hyperbolic. Dehn functions of
nilpotent groups are bounded by a polynomial \cite{Short}. Another
source of groups with polynomial Dehn function is the class of
groups with simply connected asymptotic cones (Gromov,
\cite{GromovAs}). Moreover if the asymptotic cones are simply
connected then the isodiametric function of the group is linear.

Recall that the asymptotic cone of a group $G$ (see \cite{GromovAs},
\cite{Dr} or \cite{DS}) is the ultra-limit (or Gromov-Hausdorff
limit) of a sequence $X/d_i$ where $X$ is a Cayley graph of $G$,
$\lim d_i=\infty$, $X/d_i$ is the metric space $X$ with distance
function divided by $d_i$ \cite{G2}. Asymptotic cones capture
``global" geometric properties of the group $G$.

Of a particular interest are groups with quadratic Dehn function.
That class of groups includes the classes of automatic groups and
CAT(0)-groups. Higher dimensional Heisenberg groups \cite{Alc},
\cite{OSheis} and some solvable non-virtually nilpotent groups
\cite{Dr2} also have quadratic Dehn functions. That class contains
more complicated groups as well. The most striking example so far is
the R. Thompson group $$F=\la x_0, x_1\mid x_1^{x_0^2}=x_1^{x_0x_1},
x_1^{x_0^3}=x_1^{x_0^2x_1}\ra$$ where $a^b=b\iv ab$. Recall that $F$
is the group of all piecewise linear increasing self-homeomorphisms
of the unit interval with finitely many dyadic singular points and
all slopes - powers of $2$. Guba showed in \cite{Guba} that $F$ has
a quadratic Dehn function. One of the most interesting unsolved
problems about this class is whether $\mathrm{SL}_n(\Z)$ belongs to
it for $n\ge 4$.

A very non-trivial result of Bridson and Groves \cite{BG} shows that
every cyclic extension of a finitely generated free group has
quadratic Dehn function. On the other hand, Olshanskii and I proved
\cite{OSsmalldehn} that HNN extensions of free groups having
undecidable conjugacy problem must have Dehn function at least
$n^2\log n$. Together with the result of Bridson and Groves it gives
another proof of decidability of the conjugacy problem for cyclic
extensions of free groups \cite{BMMV}.

It is still unknown whether every group with quadratic Dehn function
has decidable conjugacy problem. Olshanskii and I gave a
``quasi-proof" of that in \cite{OSsmalldehn}.

I think that it is most probable that the class of Dehn functions
$\ge n^2\log n$ is as wide as the class of time functions of Turing
machines. The next theorem confirms that conjecture in the case of
Dehn functions $\ge n^4$.

\begin{theorem} (See \cite{SBR}) \label{1}
1. Let ${\cal D}_4$ be the set of all Dehn functions $d(n)\ge n^4$
of finitely presented groups. Let ${\cal T}_4$ be the set of time
functions $t(n)\geq n^4$ of arbitrary Turing machines. Let ${\cal
T}^4$ be the set of superadditive functions which are fourth powers
of time functions. Then $${\cal T}^4\subseteq {\cal D}_4\subseteq
{\cal T}_4.$$

2. For every time function $T(n)$ of a non-deterministic Turing
machine with superadditive $T^4(n)$ there exists a finitely
presented group $G$ with Dehn function $T^4(n)$ and the isodiametric
function $T^3(n)$.
\end{theorem}

Recall that a function $f$ is superadditive if $f(n+m)\ge f(m)+f(n)$
for any $m,n$. The question of whether all Dehn functions are
superadditive is one of the unsolved mysteries of the subject.
Together with Victor Guba \cite{GS}, we proved that {\em the Dehn
function of any non-trivial free product is superadditive.} Thus if
there are non-superadditive Dehn functions then there are groups $G$
such that $G$ and $G*\Z$ have different Dehn functions!

Theorem \ref{1} has many corollaries. For example, it implies that
the {\em isoperimetric spectrum}, i.e. the set of $\alpha$'s such
that $\lfloor n^\alpha\rfloor $ is a Dehn function, contains all
numbers $\alpha\ge 4$ whose $n$-th digit can be computed by a
deterministic Turing machine in time less than $2^{2^n}$. All
``constructible" numbers (rational numbers, algebraic numbers,
values of elementary functions at rational points, etc.) satisfy
this condition. On the other hand, Theorem \ref{2} implies that if
$\alpha$ is in the isoperimetric spectrum then the $n$-th digit of
$\alpha$ can be computed in time $\le 2^{2^{2^n}}$ (see \cite{SBR}
for details). The difference in the number of $2$'s in these
expressions, is the difference between $P$ and $NP$ in Computer
Science (if $P=NP$ then there should be two $2$'s in both
expressions).

Note that before \cite{SBR} has been submitted to Annals of
Mathematics (in 1997), only a discrete set of non-integer numbers in
the isoperimetric spectrum was known \cite{Bri}. By the time the
paper appeared in print (2002), that set increased by a dense subset
in $[2,\infty)$ \cite{BB}. Groups in \cite{Bri} with Dehn functions
$\lfloor n^\alpha\rfloor$, $\alpha\not\in\N$, have easier
presentations than groups based on $S$-machines having the same Dehn
functions, but the construction in \cite{Bri}, \cite{BB} is far from
universal, and one cannot expect anything like Theorem \ref{1}
proved using their methods.

Other applications of Theorem \ref{1} are:

\begin{itemize}
\item  the first example of a finitely presented group with
NP-complete word problem,
\item examples of finitely presented
groups with easy word problem (solvable in quadratic time) and
arbitrary large (recursive) Dehn functions.
\end{itemize}

\subsection{The proof} \label{proof}

\smallskip

Here is how our construction from \cite{SBR} works. Take any Turing
machine $M$. Let $M'$ be the symmetric Turing machine described
above. Let $S(M')$ be the $S$-machine obtained as a composition of
$S'(M')$ with a positivity checking $S$-machine from \cite{SBR}
working in quadratic time. The time function of $S(M')$ is $T^3$ and
the space function is $T$ where $T$ is the time function of $M$. We
can assume that the accepting configuration of $S(M')$ is some fixed
word $W$ of the form $k_1w_1k_2w_2...k_N$ where $N>8$ (for some
small cancellation reasons) and all $w_i$ are copies of each other
written in disjoint alphabets and containing no tape letters. That
can be achieved by taking $N$ copies of the initial Turing machine
and making all of them work in parallel. Finally add one {\em hub}
relation $W=1$ to the $S$-machine $S(M')$. The resulting group
$G(M)$ has Dehn function $T^4$ provided $T^4$ is superadditive, and
isodiametric function $T^3$.

The main idea of the proof is the following. Take the standard
trapezium corresponding to a computation $W_1\to...\to W_n=W$,
identify its left and right sides (which have the same label). The
resulting diagram has one hole with boundary label $W$. Insert the
cell corresponding to the hub relation $W=1$ into the hole. The
result is a \vk diagram, called a {\em disc corresponding to the
equality $W_1=1$}. The area of that diagram is equal to the area of
the trapezium (plus 1). So it is equivalent to the product of the
time of the computation by its space. The diameter of the disc with
perimeter $\le n$ is the time of the computation. Hence the worst
area we can get is $T^4$, and the worst diameter is $T^3$. That
gives the lower bound of the Dehn function and the isodiametric
function. The upper bound is obtained by using certain surgeries on
\vk diagrams. It turns out that  every \vk diagram over the group
$G(M)$ can be decomposed into a few discs and a diagram whose area
is at most cubic (with respect to the perimeter of the original
diagram). Thus if the area of a \vk diagram is large then most of
the area is concentrated in the discs. It turns out also that the
sum of the perimeters of the discs does not exceed a constant
multiple of the perimeter of the diagram. This gives the desired
upper bound of $T^4$ for the Dehn function (it is in this part of
the proof where the superadditivity of $T^4$ is used) and $T^3$ for
the isodiametric function.

\subsection{The Dehn functions of $S$-machines and chord diagrams}

It is easy to see that the Dehn function of an $S$-machine is at
most cubic. Indeed, every \vk diagram with perimeter of length $n$
over the presentation of an $S$-machine is covered by $\theta$-bands
that start and end on the boundary. There are also $q$-bands
composed of $(q,\theta)$-cells, and $a$-bands composed of the
commutativity $(a,\theta)$-cells. It can be proved that every
$\theta$-band intersects a $q$-band (an $a$-band) at most once.
Hence the total number of $(q,\theta)$-cells is at most $n^2$. Every
other cell is an $(a,\theta)$-cell. Each $a$-band starts on the
boundary of the diagram or on the boundary of a $(q,\theta)$-cell,
so the total number of such bands is at most $n^2$ and the length of
each of them is at most $n$. Hence the total area is at most $n^3$.

It was conjectured by Rips and myself that the Dehn function of an
$S$-machine should in fact depend on the program of the $S$-machine.
We thought that $S$-machines should provide examples of groups with
Dehn functions strictly between $n^2$ and $n^3$. It turned out to be
the case. In particular, Olshanskii and I proved in
\cite{OSsmalldehn} that if $\sss$ is any $S$-machine accepting
language $L$, then the composition $\sss\circ Z(A)$ of $\sss$ and
the adding machine has Dehn function at most $n^2\log n$ and accepts
the same language $L$. Hence we get the following result.

\begin{theorem}[\cite{OSsmalldehn}]\label{3} There exists an $S$-machine with
undecidable conjugacy problem and Dehn function $n^2\log n$.
\end{theorem}

The idea of analyzing \vk diagrams over $S$-machines is to show that
if the area of a diagram is large then ``most" of the area is inside
large subtrapezia, and then analyze trapezia (i.e. computations of
the $S$-machines). Note that in a trapezium, every $\theta$-band
intersects every $q$-band, thus the band structure of a trapezium is
somewhat regular. In order to analyze irregular diagrams, Olshanskii
introduced a measure of irregularity, the {\em dispersion}. In fact,
the dispersion is an invariant of the {\em cord diagram} associated
with every \vk diagram over an $S$-machine: the role of chords is
played by the $\theta$-bands ($T$-chords) and the $q$-bands
($Q$-chords). It is similar to a Vassiliev invariant of knots.

As I mentioned before, $n^2\log n$ is the smallest Dehn function of
an HNN extension of a free group with undecidable conjugacy problem.
If the undecidability condition is dropped, one can construct Dehn
functions strictly between $n^2$ and $n^2\log n$. In particular,
Olshanskii \cite{Ol} constructed an {\em $S$-machine with
non-quadratic Dehn function bounded from above by a quadratic
function on arbitrary long intervals.} This gives the first example
of a finitely presented group with two non-homeomorphic asymptotic
cones \cite{OScones}: {\em one of the asymptotic cones of this group
is simply connected, and another one - is not.}

Non-finitely presented groups with ``very many" asymptotic cones are
constructed in \cite{DS} using completely different methods.

\begin{theorem}[Dru\c tu, Sapir \cite{DS}] \label{48} There exist finitely
generated groups with continuously many (maximal theoretically
possible if the Continuum Hypothesis is true \cite{KSTT})
non-homeomorphic asymptotic cones.
\end{theorem}

It is very interesting whether one can replace ``finitely generated"
by ``finitely presented" in Theorem \ref{48}. One can try to use
Higman embeddings from Section \ref{Hig} to construct such examples.

\subsection{Non-simply connected asymptotic cones}

Note that Theorem \ref{1} gave some of the first examples of groups
with polynomial Dehn function and non-simply connected asymptotic
cones because their Dehn functions can be polynomial (if the
original Turing machine had polynomial time function) while their
isodiametric functions are not linear. The first examples of groups
with polynomial (cubic) Dehn functions, linear isodiametric
functions and non-simply connected asymptotic cones were given in
\cite{OSnonsc}. That answered a question of Dru\c tu from \cite{Dr}.

The groups in \cite{OSnonsc} are $S$-machines. The easiest example
is this: $$G = \la \theta_1, \theta_2, a, k \mid  a^{\theta_i} = a,
k^{\theta_i} = ka, i = 1, 2\ra.$$ The $S$-machine has one tape
letter, one state letter and two rules $[k\to ka]$ (and their
inverses).

There is also an $S$-machine with Dehn function $n^2\log n$
satisfying the same asymptotic properties. Note that $n^2\log n$
cannot be lowered to $n^2$ because of a result of Papasoglu
\cite{Pap}: {\em all groups with quadratic Dehn functions have
simply connected asymptotic cones.}

If a group has non-simply connected asymptotic cone, it is natural
to ask what is its fundamental group. We do not know what are the
fundamental groups of asymptotic cones of $S$-machines. These groups
may provide some interesting invariants of $S$-machines and Turing
machines, so it is worthwhile studying them.

The following theorem gives a partial answer to the question of what
kind of groups can be fundamental groups of asymptotic cones of
finitely generated groups.

\begin{theorem} [\cite{DS}] For every
countable group $C$ there exists an
asymptotic cone of a finitely generated group $G$ whose fundamental
group is isomorphic to the free product of continuously many copies
of $C$.
\end{theorem}

The proof does not use $S$-machines but uses some small cancellation
arguments. It would be interesting to find finitely presented groups
with similar ``arbitrary" fundamental groups of asymptotic cones.
Perhaps the Higman embeddings discussed in the next section will
help solving that problem. Another very interesting problem (due to
Gromov \cite{GromovAs}) is whether there exists an asymptotic cone
of a finitely generated group with non-trivial but at most countable
fundamental group.

\section{Higman embeddings}
\label{Hig}

 The flexibility of $S$-machines allowed us to construct
several versions of Higman embeddings (embeddings of recursively
presented groups into finitely presented ones) preserving certain
properties of the group.

\subsection{An easy construction} \label{ec} The easiest known construction
of a Higman embedding is the following. Let $H$ be a recursively
presented group $\la X\mid R\ra$. Then the set of all words in
$X\cup X\iv$ that are equal to 1 in $H$ is recursively enumerable.
Hence we can assume that $R$ consists of all these words. Then there
exists an $S$-machine recognizing $R$. More precisely, for every
word $w$ in $X$, it starts with a word $q_1wq_2q_3...q_m$ and ends
with a word $\bar q_1\bar q_2...\bar q_m$ if and only if $w\in R$.

Again, as in the proof of Theorem \ref{1}, we consider $N>8$ copies
of this $S$-machine and assume that the input of the $S$-machine
$\sss$ has the form
$$K(w)=k_1q_1wq_2...q_mk_2q_1'w'q'_2...q_m'k_3...k_{N+1}$$ and the
accepting configuration $$W=k_1\bar q_1\bar q_2...\bar q_mk_2\bar
q_1'...k_{N+1}.$$ Here $w', w'',...$ are copies of $w$ written in
disjoint alphabets.

Let $G$ be the group constructed as in the proof of Theorem \ref{1}
(see Section \ref{proof}) by imposing the hub relation $W=1$ on
$\sss$. Then the word $K(w)$ is equal to 1 in $G$ if and only if
$w\in R$.

Consider now another $S$-machine $\sss'$ with input configuration
$$K'(w)=k_1q_1q_2...q_mk_2q_1'w'q'_2...q_m'k_3...k_{N+1}.$$
That machine works exactly like $\sss$ in the part of the word
between $k_2$ and $k_{N+1}$, and does nothing in the part between
$k_1$ and $k_2$. Let $G'$ be the group obtained by imposing the
relation $W=1$ on $S'$. Then $K'(w)=1$ in $G'$ if and only if $w\in
R$.

Finally consider the amalgamated product ${\cal G}=G*_AG'$ where $A$
is generated by all state and tape letters that appear in $K'(W)$.
In that group, for every $w\in R$, both $K(w)=1$ and $K'(w)=1$.
Hence $w=1$. Thus there exists a natural homomorphism from $H$ into
${\cal G}$. It is possible (and not too hard) to prove that this
homomorphism is injective. Hence $H$ is inside a finitely presented
group $G$.

Another version of embedding used in \cite{BORS} employs the so
called {\em Aanderaa trick} \cite{Aa}: instead of the amalgamated
product, we used an HNN extension (see also the survey
\cite{OStalk}).

Ones the embedding is established, it is important to understand
which properties of a group $H$ can be preserved.

\subsection{Dehn functions and quasi-isometric Higman embeddings}

First results have been obtained by Clapham \cite{Cla} and Valiev
\cite{Va} (see \cite{OScol} for the history of these results): they
proved that the solvability (even r.e. degree) of the word problem
and the level in the polynomial hierarchy of the word problem is
preserved under some versions of Higman embedding.

In \cite{BORS}, Birget, Olshanskii, Rips and the author of this
paper obtained a much stronger result.

\begin{theorem}[\cite{BORS}]
Let $H$ be a finitely generated group with word problem solvable by
a non-deterministic Turing machine with time function $\le T(n)$
such that $T(n)^4$ is superadditive. Then $H$ can be embedded into a
finitely presented group $G$ with Dehn function  $\le n^2T(n^2)^4$
in such a way that $H$ has bounded distortion in $G$.  \label{th1}
\end{theorem}

This theorem immediately implies the following characterization of
groups with word problem in $NP$.

\begin{theorem}[\cite{BORS}]\label{6}
A finitely generated group $H$ has word problem in NP if and only if
$H$ is embedded quasi-isometrically into a finitely presented group
with polynomial Dehn function.
\end{theorem}

Note that the ``if" part of this theorem is trivial: if a finitely
generated group is a (not necessarily quasi-isometric) subgroup of a
group with polynomial Dehn function, its word problem is in NP. The
converse part is highly non-trivial, although one can prove that the
embedding described in Section \ref{ec} satisfies the desired
properties (in \cite{BORS}, we used the Aanderaa trick).

From the logic point of view, Theorem \ref{6} means that for every
(arbitrary clever) algorithm solving the word problem in a finitely
generated group, there exists a finitely presented group $G>H$ such
that the word problem in $H$ (and, moreover, in $G$) can be solved
by the Miller machine $M(G)$ in approximately the same time as the
initial algorithm.

\subsection{Preserving the solvability of the conjugacy problem}

The conjugacy problem turned out to be much harder to preserve under
embeddings. Collins and Miller \cite{CM} and Gorjaga and
Kirkinski\u\i  \cite{GK} proved that even subgroups of index 2 of
finitely presented groups do not inherit solvability or
unsolvability of the conjugacy problem.

In 1976 D. Collins \cite{KT} posed the following question (problem
5.22): {\em Does there exist a version of the Higman embedding
theorem in which the degree of unsolvability of the conjugacy
problem is preserved?} In \cite{OScol}, \cite{OScol1} we solved this
problem affirmatively. In particular, we proved the following
results.

\begin{theorem}[\cite{OScol}]\label{7}
A finitely generated group $H$ has solvable conjugacy problem if and
only if it is Frattini embedded into a finitely presented group $G$
with solvable conjugacy problem.
\end{theorem}

\begin{theorem}[\cite{OScol1}]\label{8} Every countable recursively presented
group with solvable word and power problems is embeddable into a
finitely presented group with solvable conjugacy and power problem.
\end{theorem}

Recall that a subgroup $H$ of a group $G$ is Frattini embedded in
$G$ if every two elements of $H$ that are conjugate in $G$ are also
conjugate inside $H$. We say that $G$ has solvable {\em power
problem} if there exists an algorithm which, given $u, v$ in $G$
says if $v=u^n$ for some $n\ne 0$.

Theorem \ref{8} is a relatively easy application of Theorem \ref{7}.

The construction in \cite{OScol} is much more complicated than in
\cite{BORS}. First we embed $H$ into a finitely presented group
$H_1$ preserving the solvability of the word problem. Then we use
the Miller $S$-machine $M(H_1)$ to solve the word problem in $H$. In
order to overcome technical difficulties, we needed certain parts of
words appearing the computation to be always positive. The standard
positivity checkers do not work because they are $S$-machines as
well, and can insert negative letters! So we used some ideas from
the original Boone-Novikov proofs. That required introducing new
generators, $x$-letters (in addition to the $a$-, $q$-, and
$\theta$-letters in $S$-machines) and Baumslag-Solitar relations. In
addition, to analyze the conjugacy problem in $G$, we had to
consider annular diagrams which are more complicated than \vk disc
diagrams. Different types of annular diagrams (spirals, roles, etc.)
required different treatment.

We do not have any reduction of the complexity of the conjugacy
problem in $H$ to the complexity of the conjugacy problem in $G$. In
particular, solving the conjugacy problem in $G$, in some cases
required solving systems of equations in free groups (i.e. the
Makanin-Razborov algorithm).

\section{Non-amenable finitely presented groups}

One of the most important applications of $S$-machines and Higman
embeddings so far was the construction of a finitely presented
counterexample to the von Neumann problem, i.e. a finitely presented
non-amenable group without non-Abelian free subgroups \cite{OSamen}.

\subsection{Short history of the problem} Hausdorff \cite{Haus}
proved in 1914 that one can subdivide the 2-sphere minus a countable
set of points into 3 parts $A$, $B$, $C$, such that each of these
three parts can be obtained from each of the other two parts by a
rotation, and the union of two of these parts can be obtained by
rotating the third part. This implied that one cannot define a
finitely additive measure on the 2-sphere which is invariant under
the group $SO(3)$. In 1924 Banach and Tarski \cite{BT} generalized
Hausdorff's result by proving, in particular, that in
$\mathbb{R}^3$, every two bounded sets $A, B$ with non-empty
interiors can be decomposed $A=\bigcup_{i=1}^n A_i$,
$B=\bigcup_{i=1}^n B_i$ such that $A_i$ can be rotated to $B_i$,
$i=1,...,n$ (the so called Banach-Tarski paradox). Von Neumann
\cite{vN} was first who noticed that the cause of the Banach-Tarski
paradox is not the geometry of $\mathbb{R}^3$ but an algebraic
property of the group $\mathrm{SO}(3)$. He introduced the concept of
an amenable group (he called such groups ``measurable") as a group
$G$ which has a left invariant finitely additive measure $\mu$,
$\mu(G)=1$, noticed that if a group is amenable then any set it acts
upon freely also has an invariant measure and proved that a group is
not amenable provided it contains a free non-Abelian subgroup. He
also showed that groups like $\mathrm{PSL}(2, \mathbb{Z})$,
$\mathrm{SL}(2,\mathbb{Z})$ contain free non-Abelian subgroups. So
analogs of Banach-Tarski paradox can be found in $\mathbb{R}^2$ and
even $\mathbb{R}$ (for a suitable group of ``symmetries"). Von
Neumann showed that the class of amenable groups contains Abelian
groups, finite groups and is closed under taking subgroups,
extensions, and infinite unions of increasing sequences of groups.
Day \cite{day} and Specht \cite{Specht} showed that this class is
closed under homomorphic images. The class of groups without free
non-Abelian subgroups  is also closed under these operations and
contains Abelian and finite groups.

The problem of existence of non-amenable groups without non-Abelian
free subgroups probably goes back to von Neumann and became known as
the ``von Neumann problem" in the fifties. Probably the first paper
where this problem was formulated was the paper by Day \cite{day}.
It is also mentioned in the monograph by Greenleaf \cite{Greenleaf}
based on his lectures given in Berkeley in 1967. Tits \cite{Tits}
proved that every non-amenable matrix group over a field of
characteristic $0$ contains a non-Abelian free subgroup. In
particular every semisimple Lie group over a field of characteristic
$0$ contains such a subgroup.

First counterexamples to the von Neumann problem were constructed by
Olshanskii \cite{OlAmen}. He proved that the Tarsky monsters, both
torsion-free and torsion (see \cite{book}), are not amenable. Later
Adian \cite{Adian} showed that the non-cyclic free Burnside group of
odd exponent $n\ge665$ with at least two generators (that is the
group given by the presentation $\langle a_1,...,a_m\ |\ u^n=1,$
 where $u$ runs over all words in the alphabet
$\{a_1,...,a_m\} \rangle$) is not amenable.

Both Olshanskii's and Adian's examples are not finitely presented:
in the modern terminology these groups are inductive limits of word
hyperbolic groups, but they are not hyperbolic themselves. Since
many mathematicians are mostly interested in groups acting ``nicely"
on manifolds, it is natural to ask if there exists a finitely
presented non-amenable group without non-Abelian free subgroups.
This question was explicitly formulated, for example, by Grigorchuk
in \cite{KT} and by Cohen in \cite{Cohen}. This question is one of a
series of similar questions about finding finitely presented
``monsters", i.e. groups with unusual properties. Probably the most
famous problem in that series is the (still open) problem about
finding a finitely presented infinite torsion group. Other similar
problems ask for finitely presented divisible group (group where
every element has roots of every degree), finitely presented Tarski
monster, etc. In each case a finitely generated example can be
constructed as a limit of hyperbolic groups (see \cite{book}), and
there is no hope to construct finitely presented examples as such
limits.

One difficulty in constructing a finitely presented non-amenable
group without free non-Abelian subgroups is that there are ``very
few" known finitely presented groups without free non-Abelian
subgroups. Most non-trivial examples are solvable or ``almost"
solvable (see \cite{KhSap}), and so they are amenable. The only
previously known example of a finitely presented group without free
non-Abelian subgroups for which the problem of amenability is
non-trivial, is R.Thompson's group $F$ (for the definition of $F$
look in Section \ref{td}). The question of whether $F$ is not
amenable was formulated by R. Geoghegan in 1979. A considerable
amount of work has been done to answer this question but it is still
open.

\subsection{The result}

Together with A. Olshanskii, we proved the following theorem.

\begin{theorem}[\cite{OSamen}] \label{main} For every sufficiently
large odd $n$, there exists a finitely presented group \label{g1}
${\cal G}$ which satisfies the following conditions.
\begin{enumerate}
\item\label{22} ${\cal G}$ is an ascending HNN extension of a finitely generated
infinite group of exponent $n$.
\item\label{33} ${\cal G}$ is an extension of a non-locally finite group of
exponent $n$ by an infinite cyclic group.
\item\label{44} ${\cal G}$ contains a subgroup isomorphic to a free
Burnside group of exponent $n$ with $2$ generators.
\item\label{55} ${\cal G}$ is a non-amenable finitely presented group
without free non-cyclic subgroups.
\end{enumerate}
\end{theorem}

Notice that parts \ref{22} and \ref{44} of Theorem \ref{main}
immediately imply part \ref{33}. By a theorem of Adian \cite{Adian},
part \ref{44} implies that ${\cal G}$ is not amenable. Thus parts
\ref{22} and \ref{44} imply part \ref{55}.

Note that the first example of a finitely presented group which is a
cyclic extension of an infinite torsion group was constructed by
Grigorchuk \cite{Grig}. But the torsion subgroup in Grigorchuk's
group does not have a bounded exponent and his group is amenable (it
was the first example of a finitely presented amenable but not
elementary amenable group).

\subsection{The proof}

Let us present the main ideas of our construction. We first embed
the free Burnside group \label{bmn}$B(m,n)=\langle {\cal B}\rangle$
of odd exponent $n>>1$ with $m>1$ generators $\{b_1,...,b_m\}={\cal
B}$ into a finitely presented group ${\cal G}'=\langle {\cal C}\mid
{\cal R}\rangle$ where ${\cal B}\subset {\cal C}$. This is done as
in Section \ref{ec} using an $S$-machine recognizing all words of
the form $u^n$. The advantage of $S$-machines is that such an
$S$-machine can be easily and explicitly constructed (see
\cite{OStalk}). Then we take a copy \label{a11}${\cal
A}=\{a_1,...,a_m\}$ of the set ${\cal B}$, and a new generator ${\bf
t}$, and consider the group given by generators ${\cal C}\cup {\cal
A}$ and the following three sets of relations.

\begin{enumerate}
\item[(1)] \label{rrel} the set ${\cal R}$ of the relations of
the finitely presented group ${\cal G}'$ containing $B(m,n)$;
\item[(2)] \label{urel}($u$-relations) $y=u_y$, where $u_y, y\in {\cal C},$
is a certain word in ${\cal A}$ these words satisfy a very strong
small cancellation condition; these relations make ${\cal G}'$ (and
$B(m,n)$) embedded into a finitely presented group generated by
${\cal A}$;
\item[(3)] \label{rhorel} (${\bf t}$-relations) ${\bf t}\iv a_i{\bf t} =b_i,
i=1,...,m$; these relations make $\la {\cal A}\ra$ a conjugate of
its subgroup of exponent $n$ (of course, the group $\la {\cal A}\ra$
gets factorized).
\end{enumerate}
The resulting group ${\cal G}$ is obviously generated by the set
${\cal A}\cup \{{\bf t}\}$ and is an ascending HNN extension of its
subgroup $\langle {\cal A}\rangle$ with the stable letter ${\bf t}$.
Every element in $\la {\cal A}\ra$ is a conjugate of an element of
$\langle {\cal B}\rangle$, so $\la {\cal A}\ra$ is an $m$-generated
group of exponent $n$. This immediately implies that ${\cal G}$  is
an extension of a group of exponent $n$ (the union of increasing
sequence of subgroups ${\bf t}^s\la {\cal A}\ra {\bf t}^{-s},
s=1,2,...)$ by a cyclic group.

Hence it remains to prove that $\la {\cal A}\ra$ contains a copy of
the free Burnside group $B(2,n)$.

In order to prove that, we construct a list of defining relations of
the subgroup $\la {\cal A}\ra$. As we have pointed out, the subgroup
$\la {\cal A}\cup {\cal C}\ra=\la {\cal A}\ra$ of ${\cal G}$ clearly
satisfies all  {\em Burnside relations} of the form $v^n=1$. Thus we
can add all Burnside relations.

\begin{enumerate}
\item[(4)] \label{berrel} $v^n=1$ where $v$ is a word in ${\cal A}\cup {\cal C}$.
\end{enumerate}
to the presentation of group ${\cal G}$ without changing the group.

If Burnside relations were the only relations in ${\cal G}$ among
letters from ${\cal B}$, the subgroup of ${\cal G}$ generated by
${\cal B}$ would be isomorphic to the free Burnside group $B(m,n)$
and that would be the end of the story. Unfortunately there are many
more relations in the subgroup $\la {\cal B}\ra$ of ${\cal G}$.
Indeed, take any relation $r(y_1,...,y_s)$, $y_i\in {\cal C}$, of
${\cal G}$. Using $u$-relations (2), we can rewrite it as
$r(u_1,...,u_s)=1$ where $u_i\equiv u_{y_i}$. Then using ${\bf
t}$-relations, we can substitute each letter $a_j$ in each $u_i$ by
the corresponding letter $b_j\in {\cal B}$. This gives us a relation
$r'=1$ which will be called a relation {\em derived} from the
relation $r=1$, the operator producing derived relations will be
called the {\em ${\bf t}$-operator}. We can apply the ${\bf
t}$-operator again and again producing the second, third, ...,
derivatives $r''=1,r'''=1,...$ or $r=1$. We can add all {\em derived
relations}

\begin{enumerate}
\item [(5)\label{derrel}] $r'=1, r''=1,...$ for all relations
$r\in {\cal R}$
\end{enumerate}
to the presentation of ${\cal G}$ without changing ${\cal G}$.

Now consider the group $H$ generated by ${\cal C}$ subject to the
relations (1) from ${\cal R}$, the Burnside relations (4) and the
derived relations (5). The structure of the relations of $H$
immediately implies  that $H$ contains subgroups isomorphic to
$B(2,n)$. Thus it is enough to show that the natural map from $H$ to
${\cal G}$ is an embedding.

The idea is to consider two auxiliary groups. The group ${\cal G}_1$
generated by ${\cal A}\cup {\cal C}$ subject to the relations (1)
from ${\cal R}$, $u$-relations (2), the Burnside relations (4), and
the derived relations (5). It is clear that ${\cal G}_1$ is
generated by $\cal A$ and is given by relations (1) and (5) where
every letter $y\in{\cal C}$ is replaced by the corresponding word
$u_y$ in the alphabet ${\cal A}$ plus all Burnside relations (4) in
the alphabet ${\cal A}$. Let $L$ be the normal subgroup of the free
Burnside group $B({\cal A},n)$ (freely generated by ${\cal A}$)
generated as a normal subgroup by all relators (1) from ${\cal R}$
and all derived relators (5) where letters from ${\cal C}$ are
replaced by the corresponding words $u_y$. Then ${\cal G}_1$ is
isomorphic to $B({\cal A},n)/L$.

Consider the subgroup $U$ of $B({\cal A},n)$ generated (as a
subgroup) by $\{u_y \mid y\in {\cal C}\}$. The words $u_y$, $y\in
{\cal C}$, are chosen in such a way that the subgroup $U$ is a free
Burnside group freely generated by $u_y$, $y\in {\cal C}$, and it
satisfies the \label{cepp}{\em congruence extension} property,
namely every normal subgroup of $U$ is the intersection of a normal
subgroup of $B({\cal A},n)$ with $U$.

All defining relators of ${\cal G}_1$ are inside $U$. Since $U$
satisfies the congruence extension property, the normal subgroup
$\bar L$ of $U$ generated by these relators is equal to $L\cap U$.
Hence $U/\bar L$ is a subgroup of $B({\cal A},n)/L={\cal G}_1$. But
by the choice of $U$, there exists a (natural) isomorphism between
$U$ and the free Burnside group $B({\cal C},n)$ generated by ${\cal
C}$, and this isomorphism takes $\bar L$ to the normal subgroup
generated by relators from ${\cal R}$ and the derived relations (5).
Therefore $U/\bar L$ is isomorphic to $H$ (since, by construction,
$H$ is generated by ${\cal C}$ subject to the Burnside relations,
relations from ${\cal R}$ and derived relations)! Hence $H$ is a
subgroup of ${\cal G}_1$.  Let ${\cal G}_2$ be the subgroup of $H$
generated by ${\cal B}$.

Therefore we have $${\cal G}_1 \ge H \ge {\cal G}_2.$$

Notice that the map $a_i\to b_i$, $i=1,...,m$, can be extended to a
homomorphism \label{phi1}$\phi_{1,2}: {\cal G}_1\to {\cal G}_2$.
Indeed, as we mentioned above ${\cal G}_1$ is generated by ${\cal
A}$ subject to Burnside relations, all relators from ${\cal R}$ and
all derived relators (5) where letters from ${\cal C}$ are replaced
by the corresponding words $u_y$. If we apply $\phi_{1,2}$ to these
relations, we get Burnside relations and derived relations which
hold in ${\cal G}_2 \le H$.

The main technical statement of the paper shows that $\phi_{1,2}$ is
an isomorphism, that is for every relation $w(b_1,...,b_m)$ of
${\cal G}_2$ the relation $w(a_1,...,a_m)$ holds in ${\cal G}_1$.
This implies that the HNN extension $\la {\cal G}_1,{\bf t} \mid
{\bf t}\iv {\cal G}_1{\bf t} ={\cal G}_2\ra$ is isomorphic to ${\cal
G}$. Indeed, this HNN extension is generated by ${\cal G}_1$ and
${\bf t}$, subject to relations (1), (2), (4), (5) of ${\cal G}_1$
plus relations (3). So this HNN extension is presented by relations
(1)-(5) which is the presentation of ${\cal G}$. Therefore ${\cal
G}_1$ is a subgroup of ${\cal G}$, hence $H$ is a subgroup of ${\cal
G}$ as well.

The proof of the fact that $\phi_{1,2}$ is an isomorphism requires a
detailed analysis of the group $H$. This group can be considered as
a factor-group of the group $H'$ generated by ${\cal C}$ subject to
the relations (1) from ${\cal R}$ and derived relations (5)  over
the normal subgroup generated by Burnside relations (4). In other
words, $H$ is the {\em Burnside factor} of $H'$.

Burnside factors of free groups have been studied extensively
starting with the celebrated paper by Adian and Novikov \cite{AN}.
Later Olshanskii developed a geometric method of studying these
factors in \cite{book}. These methods were extended to arbitrary
hyperbolic groups in \cite{Ol95}

 The main problem we face in this paper is that $H'$ is
``very" non-hyperbolic. In particular, the set of relations ${\cal
R}$ contains many commutativity relations, so $H'$ contains
non-cyclic torsion-free Abelian subgroups which cannot happen in a
hyperbolic group.

We use a weak form of relative hyperbolicity that does hold in $H'$.
In order to roughly explain this form of relative hyperbolicity used
in the proof, consider the following example. Let $P=F_A\times F_B$
be the direct product of two free groups of rank $m$. Then the
Burnside factor of $P$ is simply $B(m,n)\times B(m,n)$. Nevertheless
the theory of \cite{book} cannot be formally applied to $P$. Indeed,
there are arbitrarily thick rectangles corresponding to relations
$u\iv v\iv uv=1$ in the Cayley graph of $P$ so diagrams over $P$ are
not A-maps in the terminology of \cite{book} (i.e. they do not look
like hyperbolic spaces). But one can obtain the Burnside factor of
$P$ in two steps. First we factorize $F_A$ to obtain $Q=B(m,n)\times
F_B$. Since $F_A$ is free, we can simply use \cite{book} to study
this factor.

Now we consider all edges labeled by letters from $A$ in the Cayley
graph of $Q$ as $0$-edges, i.e. edges of length 0. As a result the
Cayley graph of $Q$ becomes a hyperbolic space (a tree). This allows
us to apply the theory of A-maps from \cite{book} to obtain the
Burnside factor of $Q$. In fact $Q$ is weakly relatively hyperbolic
in the sense of our paper \cite{OSamen}, i.e. it satisfies
conditions (Z1), (Z2), (Z3) from the paper. The class of groups
satisfying these conditions is very large and includes groups
corresponding to $S$-machines considered in \cite{OSamen}.

Recall that set ${\cal C}$ consists of tape letters, state letters,
and command letters. In different stages of the proof some of these
letters become $0$-letters.

Trapezia corresponding to computations of the $S$-machine play
central role in our study of the Burnside factor $H$ of $H'$. As in
\cite{book}, the main idea is to construct a graded presentation
${\cal R}'$ of the Burnside factor $H$ of $H'$ where longer
relations have higher ranks and such that every \vk diagram over the
presentation of $H'$ has the so called property A from \cite{book}.
In all diagrams over the graded presentation of $H$, cells
corresponding to the relations from ${\cal R}$ and derived relations
are considered as $0$-cells or cells of rank $1/2$, and cells
corresponding to Burnside relations from the graded presentation are
considered as cells of ranks 1, 2,.... So in these \vk diagrams
``big" Burnside cells are surrounded by ``invisible" 0-cells and
``small" cells.

The main part of property A from \cite{book} is the property that if
a diagram over ${\cal R}'$ contains two Burnside cells $\Pi_1,
\Pi_2$ connected by a rectangular {\em \ct} subdiagram $\Gamma$ of
rank 0 where the sides contained in the contours of the two Burnside
cells are ``long enough"  then these two cells cancel, that is the
union of $\Gamma$, $\Pi, \Pi'$ can be replaced by a smaller
subdiagram. This is a ``graded substitute" to the classic property
of small cancellation diagrams (where \ct subdiagrams contain no
cells).

In our case, \ct subdiagrams of rank 0 turn out to be trapezia
(after we clean them of Burnside $0$-cells), so properties of \ct
subdiagrams can be translated into properties of the machine $\sss$.

\end{document}